\documentclass[a4paper]{article}
\usepackage[english]{babel}
\usepackage[cp1251]{inputenc}
\usepackage[colorlinks,linkcolor=blue,citecolor=blue]{hyperref}
\usepackage{latexsym,amsfonts,amssymb,amsthm,amsmath,amsbsy,graphicx}
\newtheorem{theorem}{Theorem}
\newtheorem{lema}{Lemma}
\newtheorem{slid}{Corollary}
\begin{document}
\thispagestyle{myheadings} \markboth{Ukrainian Mathematical Journal,
Vol. 59, No. 4, 2007, pp. 620 - 623}{Ukrainian Mathematical Journal,
Vol. 59, No. 4, 2007, pp. 620 - 623}
\bigskip
 {\noindent \bf \large   TWO-LIMIT PROBLEMS FOR ALMOST SEMICONTINUOUS\\
PROCESSES DEFINED ON A MARKOV CHAIN}\footnotemark[1]
\footnotetext{Translated from Ukrains'kyi Matematychnyi Zhurnal,
Vol. 59, No.4, pp.555-565, April, 2007. Original article submitted
February~1, 2006. This reprint differs from the original in
pagination and typographic detail.}

\bigskip
\bigskip
 { \bf Ievgen~Karnaukh
 \footnotetext{Taras Shevchenko Kyiv National University, Kyiv.\phantom{\quad \quad }
 \href{mailto:kveugene@mail.ru}{kveugene@mail.ru}
 }}\hskip 10 cm UDC 519.21
\begin{center}
\bigskip
\begin{quotation}
\noindent  {\small We consider almost upper semi-continuous
processes defined on a finite Markov chain. The distributions of the
functionals associated with the exit from a finite interval are
studied. We also consider some modification of these processes.}
\end{quotation}\footnotetext{\bigskip \hskip 2 cm 0041-5995/07/5904-0620\;\;
\copyright\; 2007\;\; Springer Science+Business Media, Inc.}
\end{center}
\bigskip
\par Problems related to the exit of a process with independent increments from
 an interval were investigated in many works (see,
e.g.~\cite{BratiychucGusak}, \cite{Pecherskiy}). Analogous problems
were investigated for processes on a finite Markov chain under the
semicontinuity condition~\cite{Gusak2,KorolyukShurenkov}. For walks
on a countable Markov chain, a two-dimensional problem was studied
in~\cite{Bratiychuc}. In the present paper, we consider
distributions of some functionals associated with the exit from a
bounded interval for a process with independent increments on a
finite Markov chain under the assumption that this process crosses a
positive level only by exponentially distributed jumps (an almost
semicontinuous process~\cite{Gusak5}).

The distributions of overjump functionals described by integral
equations on a semi-axes are defined by the projection factorization
method using an infinitely divisible factorization (instead of
canonical factorization). In the present paper, we investigate
functionals described by integral equations on an interval that can
be extended to a semi-axes. In the solution of the extended
equation, we use the method developed by Krein in~\cite{Krein,
GohberKreyn} and probability factorization identities.

Consider a two-dimensional Markov process:
$$Z(t)=\{\xi(t),x(t)\},\quad t\ge 0,$$
\noindent where $x(t)$ is a finite irreducible nonperiodic Markov
chain with the set of states ${E}'=\{1 \ldots m\}$ and the matrix of
transition probabilities
\begin{equation*}
\mathbf{P}(t)=e^{t \mathbf{Q}},\; t\ge0,\quad
\mathbf{Q}=\mathbf{N}(\mathbf{P}-\mathbf{I}),
\end{equation*}
where $\mathbf{N}=||\delta_{kr}\,\nu_k||^{m}_{k,r=1}$, $\nu_k$ are
the parameters of exponentially distributed random variables
$\zeta_k$ (the sojourn time of $x(t)$ in the state $k$),
$\mathbf{P}=\|p_{kr}\|$ is the matrix of transition probabilities of
the imbedded chain; $\boldsymbol{\pi}=(\pi_1,\ldots,\pi_m)$ is the
stationary distribution, and $\xi(t)$ is a process with stationary
conditionally independent increments for fixed values of $x(t)$
(see~\cite[p.13]{Gusak2}).
\par The evolution of the process $Z(t)$ is described by the matrix
characteristic function:
$$\boldsymbol{ \Phi}_t(\alpha)=
\|\mathrm{E}[e^{\imath \alpha \xi(t)},x(t)=
r/x(0)=k]\|=\mathbf{E}e^{\imath \alpha \xi(t)}= e^{ t
\boldsymbol{\Psi}(\alpha)}, \boldsymbol{\Psi}(0)=\mathbf{Q}.$$
In what follows, we consider processes that have cumulant
\begin{equation}\label{eq4.1}
  \boldsymbol{\Psi }(\alpha )=\boldsymbol{\Lambda }\overline{\mathbf{F}}_0(0)
  \left(\mathbf{C}
  \left(\mathbf{C}-\imath\alpha \mathbf{I}\right)^{-1}-\mathbf{I} \right)
  +\int_{-\infty}^0\left(e^{\imath\alpha x}-1 \right)d\mathbf{K}_0(x)+\mathbf{Q},
\end{equation}
where
$$d\mathbf{K}_0(x)=\mathbf{N}d\mathbf{F}(x)+\boldsymbol{\Pi}(dx),\;
\mathbf{F}(x)=\|\mathrm{P}\{\chi_{kr}<x\, ; x(\zeta _1
)=r/x(0)=k\}\|,$$ $\chi _{kr}$ are the jumps of $\xi (t)$ at the
time of transition of $x(t)$ from the state $k$ to the state $r$,
$$\boldsymbol{\Pi}(dx)=\boldsymbol{\Lambda }d\mathbf{F}_0(x),\;
\mathbf{F}_0(x)=\|\delta _{kr}F^0_k(x)\|,$$
 $F^0_k(x)$ are the distribution functions of the jumps of $\xi (t)$ if
 $x(t)=k$,
$\boldsymbol{\Lambda}=\|\delta _{kr}\lambda _k\|$, $\lambda _k$ are
the parameters of exponentially distributed random variables
$\zeta'_k$ (the time interval between two neighboring jumps of $\xi
(t) $ if $x(t)=k$), $\mathbf{C}=\left |\left |\delta
_{kr}c_{k}\right |\right |$, and $c_k$ are the parameters of
exponentially distributed positive jumps of  $\xi (t)$ if $x(t)=k$.
The process $Z(t)$ with this cumulant is the almost
upper-semicontinuous process defined in~\cite[p.43]{Gusak5}.
\par Let $\theta _s$ denote an exponentially distributed random variable with
parameter $s>0$ $(\mathrm{P}\{\theta _s>t\}=e^{-s t}, t\ge 0 )$,
independent of $Z(t)$. In this case, we rewrite the characteristic
function of $\xi (\theta _s)$ as follows
\begin{gather*}
\boldsymbol{\Phi }(s,\alpha )=\mathbf{E}e^{\imath \alpha \xi(\theta
_s)}= s\int_0^\infty e^{-st}\boldsymbol{\Phi}_t(\alpha
)dt=s\left(s\mathbf{I}-\boldsymbol{ \Psi }(\alpha )\right)^{-1},
\\
\mathbf{P}_s=s\int_0^{\infty}
e^{-st}\mathbf{P}(t)dt=\boldsymbol{\Phi}(s,0)=
s\left(s\mathbf{I}-\mathbf{Q}\right)^{-1}.
\end{gather*}
\par Denote the time of the first hit of a positive (negative) level by
$$ \tau^+(x)=\inf\{t>0:\xi (t)>x\}, x>0,$$
$$(\tau^-(x)= \inf\{t>0:\xi (t)<x\}, x<0)$$
and the time of the first exit from the interval $\left (
x-T,x\right )$, $0<x<T$, $T>0$ by:
$$\tau (x,T)=\inf\left\{t>0: \xi (t)\notin (x-T,x)\right\}.$$
We introduce the events
$$A_+(x)=\left\{\omega : \xi (\tau (x,T))\geq x\right\},
\;A_-(x)=\left\{\omega : \xi (\tau (x,T))\leq x-T\right\}.$$ Then,
for $x>0$, we can write:
\begin{equation*}
  \tau (x,T)\doteq\begin{cases}
  \tau ^+(x,T)=\tau ^+(x),\quad \omega \in A_+(x);\\
  \tau ^-(x,T)=\tau ^-(x-T),\;\omega \in A_-(x).
  \end{cases}
\end{equation*}
Denote the overjumps at the time of exit from the interval by:
\begin{gather*}
\gamma ^-_T(x)=x-T-\xi (\tau^-(x,T)),\; \gamma ^+_T(x)=\xi (\tau ^+(x,T))-x.
\end{gather*}
\par In the first part of the paper, we obtain closed-form representation
 of the following moment generating functions:
\begin{gather*}
\begin{split}
\mathbf{B}^T(s,x)&=\left \| \mathrm{E}\,\left [  e^{-s\tau (x,T)},\,\xi (\tau (x,T))\geq x,
 \,x(\tau (x,T))=r/x(0)=k\right ] \right \|=
\\&=\mathbf{E}\,\left [  e^{-s\tau ^+(x,T)},\,A_+(x)\right ],\\
\mathbf{B}_T(s,x)&=\left \| \mathrm{E}\,\left [  e^{-s\tau (x,T)},
\,\xi (\tau (x,T))\leq  x-T,\,x(\tau (x,T))=r/x(0)=k\right ] \right \|=
\\&=\mathbf{E}\,\left [  e^{-s\tau ^-(x,T)},\,A_-(x)\right ],
\end{split}\\
\mathbf{B}(s,x,T)=\mathbf{E}\,e^{-s \tau (x,T)} ,\quad
\mathbf{V}(s,\alpha ,x,T)=\mathbf{E}\,\left [ e^{i\alpha \xi
(\theta_s)},\,\tau (x,T)>\theta_s\right ],
\end{gather*}
\begin{gather*}
\mathbf{V}^{\pm}(s,\alpha ,x,T)=\mathbf{E}\,\left [  e^{i\alpha
\gamma ^{\pm}_T(x)-s\tau ^{\pm}(x,T)},\,A_{\pm}(x)\right ],
\\
\mathbf{V}_{\pm}(s,\alpha ,x,T)=\mathbf{E}\,\left [  e^{i\alpha
\xi (\tau ^{\pm }(x,T))-s\tau ^{\pm}(x,T)},\,A_{\pm}(x)\right ].
\end{gather*}
\par Denote the set of bounded functions absolutely integrable on the interval
 $I\subseteq(-\infty,\infty)$ and the set of their integral transforms by
$$\mathcal{L}_m(I)\!=\!\left\{\mathbf{G}(x)
=\|G_{kr}(x)\|\!:\!\int_{I}|G_{kr}(x)|dx<\!\infty;k,r=\overline{1,m}\right\},$$
$$\mathfrak{R}^0_m(I)=\left\{\mathbf{g}^0(\alpha )=\|g^0_{kr}(\alpha
)\|:
 g^0_{kr}(\alpha )=C_{kr}+\int_{I}e^{i\alpha x}G_{kr}(x)dx; k,r=\overline{1,m}\right\}.$$
We introduce the projection operation onto
$\mathfrak{R}^0_m((-\infty,\infty))$
\begin{gather*}
\left [  \mathbf{C}+\mathbf{g}(\alpha )\right ]_{I}=\int_{I}e^{i\alpha x}\mathbf{G}(x)dx,\;\;
\left [  \mathbf{C}+\mathbf{g}(\alpha )\right ]^0_{I}=\mathbf{C}+\int_{I}e^{i\alpha x}\mathbf{G}(x)dx,\\
\left [  \mathbf{C}+\mathbf{g}(\alpha )\right ]_{-}=\left [
         \mathbf{C}+\mathbf{g}(\alpha )\right ]_{(-\infty,0)},\;\left [  \mathbf{C}+\mathbf{g}(\alpha
)\right ]_{+}=\left [  \mathbf{C}+\mathbf{g}(\alpha )\right ]_{(0,\infty)}.
\end{gather*}
Note that $\mathbf{V}(s,\alpha ,x,T)\in \mathfrak{R}^0_m((x-T,x))$,
$\mathbf{V}^+(s,\alpha ,x,T)\in \mathfrak{R}^0_m([x,\infty))$,
$\mathbf{V}^-(s,\alpha ,x,T)\in \mathfrak{R}^0_m((-\infty,x-T])$.
\par Further, we introduce extrema of $\xi (t)$ and the
corresponding distribution functions
$$ \xi ^{\pm
}(t)=\sup\limits_{0\leq u\leq t}(\inf)\xi (u),\; \xi ^{\pm
}=\sup\limits_{0\leq u\leq \infty}(\inf)\xi (u),$$
$$\overline{\xi}(t)=\xi (t)-\xi ^+(t),\;
\stackrel{{\scriptscriptstyle\textsf{v}}}{\xi} (t)=\xi (t)-\xi
^-(t);$$
$$  \mathbf{P}_+(s,x)=\mathbf{P}\left\{\xi
  ^+(\theta_s)<x\right\},\,x>0,\;
  \mathbf{P}^-(s,x)=\mathbf{P}\left\{\overline{\xi}(\theta_s)<x\right\},\,x<0,$$
$$  \mathbf{p}_{+ }(s)=\mathbf{P}\left\{\xi ^{+ }(\theta_s)=0\right\},\;
  \mathbf{q}_{+ }(s)=\mathbf{P}_s-\mathbf{p}_{+ }(s),$$
$$  \mathbf{p}^*_+(s)=\mathbf{p}_+(s)\mathbf{P}_s^{-1},\;
  \mathbf{R}^*_+(s)=\mathbf{C}\mathbf{p}_+^*(s).$$
\begin{lema}{\rm\cite[p.49]{Gusak2}} For the two-dimensional Markov
process
 $Z(t)=\{\xi(t),x(t)\}$,  the following factorization identity is
 true
\begin{gather}
\boldsymbol{\Phi}\left(s,\alpha \right)=\mathbf{E}e^{\imath\alpha \xi (\theta _s)}= \left\{
\begin{matrix}
\boldsymbol{\Phi }_+(s,\alpha )\mathbf{P}_s^{-1}\boldsymbol{\Phi }^-(s,\alpha ),\\
\boldsymbol{ \Phi }_-(s,\alpha )\mathbf{P}_s^{-1}\boldsymbol{\Phi
}^+(s,\alpha ),
\end{matrix} \right.\label{eq4.2}
\end{gather}
where
$$
\boldsymbol{\Phi}_{\pm }(s,\alpha ) =\mathbf{E}e^{\imath\alpha
 \xi^{\pm }(\theta_s)},\;\boldsymbol{\Phi}^-(s,\alpha)=\mathbf{E}e^{\imath\alpha
 \overline{\xi}(\theta_s)},
\boldsymbol{\Phi}^+(s,\alpha)=\mathbf{E}e^{\imath\alpha
 \stackrel{{\scriptscriptstyle\textsf{v}}}{\xi}(\theta_s)}. \notag
$$
\end{lema}
\bigskip
\begin{theorem}
For a process $Z (t)$ with cumulant~\eqref{eq4.1},
$\mathbf{B}^T(s,x)$ is determined by the relation
\begin{multline}\label{eq4.3}
s\mathbf{B}^T(s,x)=s\left(\mathbf{I}-\mathbf{p}^*_+(s)\right)e^{-\mathbf{R}^*_+(s)x}
-\mathbf{p}^*_+(s)\int_{-\infty}^{x-T}d\mathbf{P}^-(s,y)e^{-\mathbf{C}(x-y)}\mathbf{C}_0^T(s)-\\
-\left(\mathbf{I}-\mathbf{p}^*_+(s)\right)\int_{0}^{x}e^{-\mathbf{R}^*_+(s)
z}
 \mathbf{R}^*_+(s)
 \int_{-\infty}^{x-z-T}d\mathbf{P}^-(s,y)e^{-\mathbf{C}(x-y-z)}dz\mathbf{C}^T_0(s),\;
 0<x<T,
\end{multline}
\begin{equation*}
\mathbf{C}_0^T (s)=\mathbf{\Lambda
}\overline{\mathbf{F}}_0(0)\left(\mathbf{I}+\mathbf{C}
\int_{0}^{T}e^{\mathbf{C}z}\mathbf{B}^T(s,z)dz \right).
\end{equation*}
\end{theorem}
\noindent \textbf{Proof.} Using the stochastic relations for $\tau _
{kr}^+(x,T)$, namely,
\begin{equation}\label{eq4.4}
  \tau ^+_{kr}(x,T)\doteq
  \begin{cases}
    \zeta'_k ,& \zeta'_k<\zeta _k, \xi_k >x, \\
    \zeta'_k +\tau ^+_{kr}(x-\xi_k ,T),&\zeta'_k<\zeta _k,x-T<\xi_k <x,\\
    \zeta _k+\tau ^+_{jr}(x-\chi _{kj},T),&\zeta'_k>\zeta _k,x-T<\chi
    _{kj}<x,
  \end{cases}
\end{equation}
where the subscripts denote, respectively, the initial state and the
state of the chain $x(t)$ at the time of exit form $(x-T,x)$
$\left(x(0)=k, x(\tau (x,T))=r \right)$, we get
\begin{equation*}
\begin{split}
B^T_{kr}(s,x)=&\lambda_k\int_{0}^{\infty}e^{-(s+\lambda_k+\nu_k)y}dy\left(
\int_{x}^{\infty}dF^0_k(z)+
\int_{x-T}^{x}dF_k^0(z)B^T_{kr}(s,x-z)\right)+\\
&+\sum_{j=1}^{m}\nu
_k\int_{0}^{\infty}e^{-(s+\lambda_k+\nu_k)y}dy\int_{x-T}^{x}dF_{kj}(z)B^T_{jr}(s,x-z),
\;0<x<T.
\end{split}
\end{equation*}
We can rewrite these equations in the matrix from:
\begin{gather}
 \left(s\mathbf{I}+\mathbf{\Lambda }+\mathbf{N} \right)\mathbf{B}^T(s,x)
 =\mathbf{\Lambda }\overline{\mathbf{F}}_0(x)+\int_{x-T}^{x}
 d\mathbf{K}_0(z)\mathbf{B}^T(s,x-z),\; 0<x<T,\label{eq4.5}\\
\mathbf{B}^T(s,x)=0,\;x\geq T,\quad \mathbf{B}^{T}(s,x)=\mathbf{I},\;x<0.\notag
\end{gather}
Performing the substitution
$\overline{\mathbf{B}}\phantom{|}^T(s,x)=\mathbf{I}-\mathbf{B}^T(s,x),$
in~\eqref{eq4.5}, we derive the following equation for
$\overline{\mathbf{B}}\phantom{|}^T(s,x)$ $\left(0<x<T \right)$
\begin{equation*}
  (s\mathbf{I}+\mathbf{\Lambda }+\mathbf{N} )\overline{\mathbf{B}}\phantom{|}^T(s,x)
  =\left(s\mathbf{I}-\mathbf{Q} \right) + \int_{-\infty}^{\infty}d\mathbf{K}_0(z)
  \overline{\mathbf{B}}\phantom{|}^T(s,x-z),\;0<x<T.
\end{equation*}
Extending this equation to the semi-axis $x>0$, we obtain
\begin{equation}\label{eq4.6}
  (s\mathbf{I}+\mathbf{\Lambda}+\mathbf{N} )\overline{\mathbf{B}}\phantom{|}^T(s,x)=
  \left(s\mathbf{I}-\mathbf{Q} \right)+\int_{-\infty}^{\infty}d\mathbf{K}_0(z)
  \overline{\mathbf{B}}\phantom{|}^T(s,x-z)+e^{-\mathbf{C}x}\mathbf{C}_0^T(s)
  I\left\{x>T\right\}.
\end{equation}
We denote ${C}_\epsilon (x)=e^{-\epsilon x}I\left\{x>0\right\}$ and
consider the following equation for $\mathbf{Y}_\epsilon (T,s,x)$ ($
x>0, \epsilon >0$), instead of Eq.~\eqref{eq4.6}:
\begin{equation}\label{eq4.7}
  (s\mathbf{I}+\mathbf{\Lambda }+\mathbf{N}  )\mathbf{Y}_\epsilon (T,s,x)=
  \left(s\mathbf{I}-\mathbf{Q} \right)C_\epsilon (x)+\int_{-\infty}^{\infty}
  \!\!d\mathbf{K}_0(z)Y_\epsilon
  (T,s,x-z)+e^{-\mathbf{C}x}\mathbf{C}_0^T(s) I\left\{x>T\right\},
\end{equation}
Applying the integral transform with respect to $x$
to~\eqref{eq4.7}, we get
\begin{equation}\label{eq4.8}
\begin{split}
  \left(s\mathbf{I}-\boldsymbol{\Psi}(\alpha ) \right)\widetilde{\mathbf{Y}}_\epsilon (T,s,\alpha )
  =&\left(s\mathbf{I}-\mathbf{Q} \right)\int_{0}^{\infty}e^{i\alpha z}
  e^{-\epsilon z}dz+\int_{0}^{\infty}e^{i\alpha z}e^{-\mathbf{C}z}\mathbf{C}_0^T(s) I\left\{z>T\right\}dz-\\
&-\left [  \widetilde{\mathbf{K}}_0(\alpha )\widetilde{\mathbf{Y}}_\epsilon (T,s,\alpha )\right ]_{-},
\end{split}
\end{equation}
$$\widetilde{\mathbf{K}}_0(\alpha )=\int_{0}^{\infty}e^{i\alpha z}d\mathbf{K}_0(z),\;
\widetilde{\mathbf{Y}}_\epsilon (T,s,\alpha
)=\int_{0}^{\infty}e^{i\alpha z}\mathbf{Y}_\epsilon (T,s,z)dz.
$$
Using~\eqref{eq4.2} and~\eqref{eq4.8} and performing the projection
$\left [ \,\right ]_{+}$, we obtain
\begin{multline*}
s \widetilde{\mathbf{Y}}_\epsilon (T,s,\alpha )=
 \boldsymbol{\Phi}_+(s,\alpha )\mathbf{P}^{-1}_s\biggl [\boldsymbol{\Phi}^-(s,\alpha )
  \biggl(\left(s\mathbf{I}-\mathbf{Q} \right)
  \int_{0}^{\infty}e^{i\alpha z}e^{-\epsilon z}dz+\\
  +\int_{0}^{\infty}e^{i\alpha z}e^{-\mathbf{C}z}
  \mathbf{C}_0^T(s) I\left\{z>T\right\}dz\biggr)  \biggr
  ]_+.
\end{multline*}
Inverting this relation, we get
\begin{equation}\label{eq4.9}
\begin{split}
  s\mathbf{Y}_\epsilon (T,s,x)=&\int_{0}^{x}d\mathbf{P}_+(s,z)\mathbf{P}^{-1}_s
  \int_{-\infty}^{0}d\mathbf{P}^{-}(s,y)\left( s\mathbf{I}-\mathbf{Q} \right)e^{-\epsilon (x-y-z)}+\\
&+\int_{0}^{x}d\mathbf{P}_+(s,z)\mathbf{P}^{-1}_s
 \int_{-\infty}^{\min\left\{x-z-T,0\right\}}d\mathbf{P}^-(s,y)e^{-\mathbf{C}(x-y-z)}\mathbf{C}_0^T(s),
\end{split}
\end{equation}
Since $\mathbf{Y}_\epsilon (T,s,x)\rightarrow
\overline{\mathbf{B}}\phantom{|}^T(s,x)$ as $\epsilon \rightarrow
0$, $0<x<T$, relation~\eqref{eq4.9} yields
$$
s\overline{\mathbf{B}}\phantom{|}^T(s,x)=\int_{0}^{x}d\mathbf{P}_+(s,z)
 \left( s\mathbf{I}-\mathbf{Q} \right)
 +\int_{0}^{x}d\mathbf{P}_+(s,z)\mathbf{P}^{-1}_s
  \int_{-\infty}^{x-z-T}d\mathbf{P}^-(s,y)e^{-\mathbf{C}(x-y-z)}\mathbf{C}^T_0(s).
$$
Taking into account that~\cite[p.45]{Gusak5}
\begin{equation}\label{eq4.11}
\mathbf{P}\left\{\xi
^+(\theta_s)>x\right\}=\left(\mathbf{I}-\mathbf{p}_+^*(s)
\right)e^{-\mathbf{R}_+^*(s)x}\mathbf{P}_s,\;x>0,
\end{equation}
we get~\eqref{eq4.3}.
\par Below, we present (without proof) analogs of the Pecherskii identities
(see~\cite[p.108]{Pecherskiy}).
\begin{lema}
The following identities are true for $Z(t)$
\begin{gather}
\mathbf{V}(s,\alpha ,x,T)=\boldsymbol{\Phi}(s,\alpha )
\left(\mathbf{I}-\mathbf{V}_+(s,\alpha ,x,T)-\mathbf{V}_-(s,\alpha ,x,T) \right),\label{eq4.12}\\
\mathbf{V}(s,\alpha ,x,T)=\boldsymbol{\Phi}_+(s,\alpha )\mathbf{P}^{-1}_s
 \left [  \boldsymbol{\Phi}^-(s,\alpha )(\mathbf{I}-\mathbf{V}_+(s,\alpha ,x,T))\right
]_{[x-T,\infty)},\label{eq4.13}\\
\mathbf{V}(s,\alpha ,x,T)=\boldsymbol{\Phi}_-(s,\alpha
)\mathbf{P}^{-1}_s \left [  \boldsymbol{\Phi}^+(s,\alpha
)(\mathbf{I}-\mathbf{V}_-(s,\alpha ,x,T))\right
]_{(-\infty,x]}.\label{eq4.14}
\end{gather}
\end{lema}
\par Denote the joint distribution of $\left\{\xi (\theta_s),\xi^+ (\theta_s),\xi^-
(\theta_s)\right\}$ by
\begin{equation*}
  \begin{split}
  \mathbf{H}_s(T,x,y)&=\left\|\mathrm{P}
   \left\{\xi (\theta_s)<y,\xi ^+(\theta _s)<x,\xi ^-(\theta_s)>x-T,\;x(\theta _s)=r/x(0)=k\right\}\right\|\\
  &=\mathbf{P}\left\{\xi (\theta_s)<y,\tau (x,T)>\theta_s\right\}.
  \end{split}
\end{equation*}
\begin{theorem}
For a process $Z (t)$ with cumulant~\eqref{eq4.1}, the joint
distributions of $\left \{\tau ^+(x,T),\gamma ^+_T(x)\right\}$ and
$\left \{\tau ^+(x,T),\xi(\tau ^+(x,T))\right\}$ are defined by the
relations
\begin{equation}\label{eq4.16}
  \begin{cases}
    \mathbf{V}^+(s,\alpha ,x,T)=\mathbf{B}^T(s,x)
    {\mathbf{C}}\left(\mathbf{C}-i\alpha\mathbf{I}  \right)^{-1},\;0<x<T, \\
\displaystyle \mathbf{V}_+(s,\alpha ,x,T)=e^{i\alpha
x}\mathbf{V}^+(s,\alpha ,x,T).
  \end{cases}
\end{equation}
The characteristic function of  $\xi (\theta_s)$ up to the time of
exit from the interval has the following form
\begin{equation}\label{eq4.17}
  \mathbf{V}(s,\alpha ,x,T)
  =\boldsymbol{\Phi}_+(s,\alpha )\mathbf{P}^{-1}_s
   \left [\boldsymbol{\Phi}^-(s,\alpha )\left(\mathbf{I}- e^{i\alpha x}
   \mathbf{B}^T(s,x)\mathbf{C}
   \left(\mathbf{C}-i\alpha \mathbf{I} \right)^{-1} \right)\right ]_{ [
   x-T,\infty)}.
\end{equation}
The corresponding distribution has density $(x-T<y<x,\;y\neq 0)$
\begin{multline}\label{eq4.18}
\mathbf{h}_s(T,x,y)=\frac{\partial}{\partial y}\mathbf{H}_s(T,x,y)=\\
= \mathbf{p}^*_+(s)(\mathbf{P}^-(s,y))'I\left\{y<0\right\}
 +\left(\mathbf{I}-\mathbf{p}^*_+(s)\right)\int_{x-T}^{\min\left\{0,y\right\}}
e^{-\mathbf{R}^*_+(s)(y-z)}\mathbf{R}^*_+(s)d\mathbf{P}^-(s,z)\\
 -\mathbf{p}^*_+(s)
  \int_{-\infty}^{y-x}d\mathbf{P}^-(s,z)\mathbf{B}^T(s,x)\mathbf{C}e^{-\mathbf{C}(y-x-z)}
-\left(\mathbf{I}-\mathbf{p}^*_+(s)\right)\int_{0}^{y-(x-T)}
 e^{-\mathbf{R}^*_+(s)v}\times\\
 \times\int_{-\infty}^{y-v-x}\mathbf{R}^*_+(s)d\mathbf{P}^-(s,z)\mathbf{B}^T(s,x)
   \mathbf{C}e^{-\mathbf{C}(y-v-x-z)}dv,
\end{multline}
with atom at zero
\begin{equation}\label{eq4.19}
\mathbf{P}\left\{\xi (\theta_s)=0,\tau (x,T)>\theta_s\right\}\!=s
 \left(s\mathbf{I}+\boldsymbol{\Lambda}-\mathbf{N}
 \left(\left \|\mathrm{P}\left\{\chi _{kr}=0,x(\zeta _1)=r/x(0)=k\right\}
  \right \|\!-\!\mathbf{I} \right) \right)\!^{-1}.
\end{equation}
The probability of nonexit from the interval $\left (x-T,x\right )$
is determined by the relation
\begin{equation}\label{eq4.20}
\mathbf{P}\left\{\tau
(x,T)>\theta_s\right\}=\int_{x-T}^{x}d\mathbf{H}_s(T,x,y).
\end{equation}
For the moment generating functions of $\tau (x,T)$ and $\tau
^-(x,T)$, the following relations are true:
\begin{gather}\label{eq4.21}
  \begin{cases}
    \mathbf{B}(s,x,T)=\mathbf{I}
    -\mathbf{P}\left\{\tau (x,T)>\theta_s\right\}\mathbf{P}^{-1}_s,\; 0<x<T,\\
    \mathbf{B}_T(s,x)=\mathbf{B}(s,x,T)-\mathbf{B}^T(s,x),\; 0<x<T.
  \end{cases}
\end{gather}
\end{theorem}
\noindent \textbf{Proof.}  Using stochastic relations~\eqref{eq4.4}
for $\tau ^+(x,T)$ and
\begin{equation*}
  \gamma  ^+_T(x)_{kr}\doteq
  \begin{cases}
    \xi _k ,& \zeta'_k<\zeta _k, \xi_k >x, \\
    \gamma ^+_T(x-\xi _k)_{kr},&\zeta'_k<\zeta _k,x-T<\xi_k <x,\\
    \gamma ^+_T(x-\chi _{kj})_{jr},&\zeta'_k>\zeta _k,x-T<\chi  _{kj}<x,
  \end{cases}
\end{equation*}
we obtain the following equation for $\gamma ^+_T(x)$:
\begin{multline*}
  \left(s\mathbf{I}+\mathbf{N}+\boldsymbol{\Lambda} \right)\mathbf{V}^+(s,\alpha ,x,T)=
  \boldsymbol{\Lambda}\overline{\mathbf{F}}_0(0)e^{-\mathbf{C}x}\mathbf{C
  }\left(\mathbf{C}-i\alpha \mathbf{I} \right)^{-1}+\int_{x-T}^{x}d\mathbf{K}_0(z)
  \mathbf{V}^+(s,\alpha ,x-z,T).
\end{multline*}
Using~\eqref{eq4.5}, we deduce the first relation in~\eqref{eq4.16}.
The second relation follows from the definition of $\gamma ^+_T(x)$.
Relation~\eqref{eq4.17} follows from~\eqref{eq4.13}. After the
inversion with respect to $\alpha$, we obtain~\eqref{eq4.18}
and~\eqref{eq4.19} from~\eqref{eq4.17}.
\par Using relations~\eqref{eq4.18},~\eqref{eq4.19} and an analog of
the Bratiichuk formulas~\cite[p.187]{BratiychucGusak}, we can obtain
a matrix analog for the generatrices of the joint distributions of
$\left\{\tau ^+(x,T),\xi (\tau ^+(x,T))\right\}$ and $\left\{\tau
^-(x,T),\xi (\tau ^-(x,T))\right\}$:
\begin{theorem}
For $Z(t)$, the following relations are true:
\begin{gather}
s\mathbf{E}\,\left [  e^{-s\tau ^+(x,T)},\xi (\tau ^+(x,T))> z,\right ]
 =\int_{x-T}^{x}d\mathbf{H}_s(T,x,y)\overline{\mathbf{K}}_0(z-y),\; z> x,\label{eq4.22}\\
s\mathbf{E}\,\left [  e^{-s\tau ^-(x,T)},\xi (\tau ^-(x,T))< z,\right ]
 =\int_{x-T}^{x}d\mathbf{H}_s(T,x,y)\mathbf{K}_0(z-y),\; z<x-T ,\label{eq4.23}\\
\overline{\mathbf{K}}_0(x)=\int_{x}^{\infty}d\mathbf{K}_0(y),x>0,\;
\mathbf{K}_0(x)=\int_{-\infty}^{x}d\mathbf{K}_0(y),x<0\notag.
\end{gather}
\end{theorem}
\noindent \textbf{Proof.} According
to~\cite[p.469]{KorolyukShurenkov}, we have
\begin{equation}\label{eq4.24}
  \mathrm{E}_i\left [  e^{{-s\tau(x,T) }}f(x-\xi (\tau(x,T) ),x(\tau(x,T) ))\right ]-f(x,i)=
  \mathrm{E}_i \int_{0}^{\tau(x,T) }e^{-st}g(x-\xi (t),x(t))dt,
\end{equation}
where $f$ is a bounded function, $g=Af-sf$ and $A$ is a generator of
the semigroup defined by the cumulant $\boldsymbol{\Psi }(\alpha )$.
For the right-hand side of the equation, we have
\begin{multline}\label{eq4.25}
\mathrm{E}_i\int_{0}^{\tau (x,T)}\!\!\!e^{-st}g(x-\xi (t),x(t))dt=
\sum_{j=1}^{m}\!\int_{0}^{\infty}\!\!e^{-st}\mathrm{E}_i\left [
g(x-\xi (t),j),\tau (x,T)>t,
    x(t)=j\right ]dt=\\
  =\sum_{j=1}^{m}s^{-1}\int_{x-T}^{x}g(x-y,j)d\left(H_s(T,x,y) \right)_{ij}.
\end{multline}
Assuming that $f(x,i)=I\{x\geq-z\}\delta_{ir}$, $z>0$, $i,r\in E'$,
we obtain
\begin{multline}\label{eq4.26}
\mathrm{E}_i\left [  e^{-s\tau (x,T)}f(x-\xi (\tau (x,T)),x(\tau (x,T)))\right ]
      -f(x,i)= \\
  =\mathrm{E}_i\left [  e^{-s\tau (x,T)},\gamma ^+_T(x)\geq z,x(\tau (x,T))=r\right ],
\end{multline}
\begin{equation}\label{eq4.27}
g(x,j)=
    \int_{-\infty}^{\infty}I\left\{x-y\leq -z\right\}dK_0^{jr}(y)=\overline{K}^{jr}_0(x+z).
\end{equation}
Substituting~\eqref{eq4.26} and~\eqref{eq4.27} into~\eqref{eq4.24}
and taking~\eqref{eq4.25} into account, we get
\begin{equation}\label{eq4.28}
\mathrm{E}_i\left [  e^{-s\tau (x,T)},\gamma ^+_T(x)\geq z,x(\tau
(x,T))=r\right ]
=\sum_{j=1}^{m}s^{-1}\int_{x-T}^{x}\overline{K}^{jr}_0(x-y+z)d\left(H_s(T,x,y)
\right)_{ij}.
\end{equation}
Using the definition of $\gamma ^+_T(x)$, we derive~\eqref{eq4.22}
from~\eqref{eq4.28}. By analogy, we obtain relation~\eqref{eq4.23}.
\par Consider the behavior of the functions $\mathbf{B}^T(s,x)$ and
$\mathbf{H}_s(T,x,y)$ as $s\rightarrow 0$. Denote $$
  \mathbf{M}(y)=\lim_{s\rightarrow 0}s^{-1}\mathbf{p}^*_+(s)\mathbf{P}^-(s,y).$$
The existence of this function follows from the reasoning presented
below. According to~\cite[p.46]{Gusak5}, for the moment generating
function of $\overline{\xi } (\theta _s)$ we have
\begin{equation*}
\lim_{s\rightarrow
0}s^{-1}{\mathbf{p}}^*_+(s)\mathbf{E}e^{r\overline{\xi }(\theta
_s)}=
 - \left({\mathbf{p}}^*_+(0)\mathbf{C}
 - r\mathbf{I}\right)\left(\mathbf{C}-r\mathbf{I}
 \right)^{-1}\boldsymbol{\Psi}^{-1}(-ir).
\end{equation*}
Then we can define $\mathbf{M}(y)$ as a function for which
\begin{equation*}
\int_{-\infty}^{0}e^{r y}d\mathbf{M}(y)=
 - \left({\mathbf{p}}^*_+(0)\mathbf{C}
 - r\mathbf{I}\right)\left(\mathbf{C}-r\mathbf{I} \right)^{-1}\boldsymbol{\Psi}^{-1}(-ir).
\end{equation*}
Furthermore, it follows from~\cite[p.41]{Gusak5} that $$
  \mathbf{p}^*_+(0)=\left(\mathbf{I}
 -\left\|\mathrm{P}\left\{\tau ^+(0)<\infty,x(\tau ^+(0))=r/x(0)=k\right\}\right\| \right).
$$
\begin{slid}
For the process $Z(t)$ following relations are true:
\begin{multline}\label{eq4.29}
\lim_{s\rightarrow 0}s^{-1}\mathbf{h}_s(T,x,y)=
\mathbf{M}'(y)I\left\{y<0\right\}+
 \left(\mathbf{I}-\mathbf{p}^*_+(0) \right)\int_{x-T}^{\min(0;y)}e^{-\mathbf{R}^*_+(0)(y-z)}
  \mathbf{C}d\mathbf{M}(z)-\\
-\left(\mathbf{I}-\mathbf{p}^*_+(0)
\right)\int_{0}^{y-(x-T)}e^{-\mathbf{R}^*_+(0)v}\mathbf{C}
  \int_{-\infty}^{y-v-x}d\mathbf{M}(z)\mathbf{B}^T(x)\mathbf{C}e^{-\mathbf{C}(y-v-x-z)}dv-\\
-\int_{-\infty}^{y-x}d\mathbf{M}(z)\mathbf{B}^T(x)\mathbf{C}e^{-\mathbf{C}(y-x-z)},
\end{multline}
\begin{multline}\label{eq4.30}
\mathbf{B}^T(x)=\lim_{s\rightarrow 0}\mathbf{B}^T(s,x) =
\left(\mathbf{I}-\mathbf{p}^*_+(0) \right)e^{-\mathbf{R}^*_+(0)x}-
   \int_{-\infty}^{x-T}d\mathbf{M}(y)e^{-\mathbf{C}(x-y)}\mathbf{C}^T_0(0)-\\
-\left(\mathbf{I}-\mathbf{p}^*_+(0)
\right)\int_{0}^{x}e^{-\mathbf{R}^*_+(0)z}\mathbf{C}
  \int_{-\infty}^{x-z-T}d\mathbf{M}(y)e^{-\mathbf{C}(x-y-z)}dz\mathbf{C}^T_0(0),
\end{multline}
$$\mathbf{C}^T_0(0)=\boldsymbol{\Lambda}\overline{\mathbf{F}}_0(0)
   \left(\mathbf{I}+\mathbf{C}\int_{0}^{T}e^{\mathbf{C}z}\mathbf{B}^T(z)dz \right).$$
\end{slid}
\par Assume that $\chi _{kr}\equiv 0$. If $x(t)=k$, $k=1,...,m$,
then we set
$$
 \xi (t)=\sum_{n\leq \nu'_k(t)}{\xi '_n} ^k -
\sum_{n\leq \nu  _k(t)}\xi^k _n,
$$
where $\nu'_k(t)$ and $\nu_k(t)$ are Poisson processes with the
rates $\lambda _k^1$ and $\lambda _k^2$, respectively. ${\xi '_n}
^k$ and $\xi_n ^k$ are independent positive random variables, ${\xi
'_n} ^k$ are exponentially distributed with the parameters $c_k$,
the variables $\xi_n ^k$ have an arbitrary distribution with bounded
expectation $m_k$. It is clear that, in this case, the process
$Z(t)=\left\{\xi (t),x(t)\right\}$ is the almost
upper-semicontinuous and has cumulant
$$
\boldsymbol{\Psi }(\alpha )=\boldsymbol{\Lambda
}\overline{\mathbf{F}}_0(0)\left(\mathbf{C}
  \left(\mathbf{C}-\imath\alpha \mathbf{I}\right)^{-1}-\mathbf{I} \right)
  +\int_{-\infty}^0\left(e^{\imath\alpha x}-\mathbf{I} \right)\boldsymbol{\Pi }(dx)
  +\mathbf{Q},
$$
where
 $$\boldsymbol{\Lambda }=\|\delta _{kr}(\lambda _k^1+\lambda
_k^2 )\|,\; \overline{\mathbf{F}}_0(0)=\|\delta _{kr}\lambda
_k^1/(\lambda _k^1+\lambda _k^2 )\|,\;
\boldsymbol{\Pi}(dx)=\boldsymbol{\Lambda
}\mathbf{F}_0(0)d\mathbf{F}^1_0(x),$$
$$\mathbf{F}_0(0)=\mathbf{I}-\overline{\mathbf{F}}_0(0),\;
\mathbf{F}^1_0(x)=\|\delta _{kr}\mathrm{P}\left\{-\xi
^k_n<x\right\}\|, x<0.$$

Consider the process $\eta _{B,u}(t)$ defined by the stochastic
relations
\begin{equation}\label{eq4.31}
  \eta ^{kr}_{B,u}(t)\doteq
  \begin{cases}
    u+\xi _{kr}(t) & t<T_1, \\
    B & t\in (T_1,T_2),T_1<\infty,\\
    \eta^{jr} _{B,B-\xi_1^j}(t-T_2)& t>T_2,x(T_2)=j,T_1<\infty,
  \end{cases}
\end{equation}
where the superscripts $kr$ mean that $x(t)=r, x(0)=k$. Note that,
$T_1\doteq\tau ^+(v)$, $v=B-u$ and $T_2\doteq T_1+\zeta _*$, $\zeta
_*$ is the time of the first negative jump independent of $T_1$.
\par The process $\eta _{B,u}(t)$ is called the risk process in a Markov environment with
stochastic premium function and bounded
reserve(see~\cite{Buhlmann}~-~\cite{Kartashov}). We also consider
the dividend process $Y_{B,u}(t)\doteq u+\xi(t)-\eta _{B,u}(t)$
(see~\cite[p.169]{Buhlmann}).
\begin{theorem}
  The distribution of  $\eta _{B,u}(\theta _s)$ is determined by the
  characteristic function
\begin{multline}\label{eq4.32}
    \boldsymbol{\Phi}_{B,u}(s,\alpha )=\mathbf{E}e^{i\alpha \eta _{B,u}(\theta_s)}=e^{i\alpha B}
    \left(
    \left(\mathbf{C}-i\alpha \mathbf{I}\right)\mathbf{p}^*_+(s)e^{-i\alpha
    v}-\right.\\
    \left.-\left(\mathbf{I}-\mathbf{p}^*_+(s) \right)
    e^{-\mathbf{R}^*_+(s)v}\mathbf{R}^*_+(s)\right)
   \left(\mathbf{R}^*_+(s)-i\alpha \mathbf{I}\right)^{-1}
   \boldsymbol{\Phi }^-(s,\alpha )+\\
    +\left(\mathbf{I}-\mathbf{p}^*_+(s) \right)
   e^{-\mathbf{R}^*_+(s)v}
   \left(s\mathbf{I}+\boldsymbol{\Lambda
    }\mathbf{F}_0(0)-\mathbf{Q} \right)^{-1}\left(se^{i\alpha B}+
    \boldsymbol{\Lambda }\mathbf{F}_0(0)
    \widetilde{\boldsymbol{\Phi }}_B(s,\alpha ) \right),
\end{multline}
\begin{multline}\label{eq4.33}
 \widetilde{\boldsymbol{\Phi }}_B(s,\alpha
)=\int_{-\infty}^{0}d\mathbf{F}^1_0(z)\boldsymbol{\Phi
}_{B,B+z}(s,\alpha ) =e^{i\alpha B}\left(\boldsymbol{\Lambda
}\mathbf{F}_0(0) \right)^{-1} \left(s\mathbf{I}+ \boldsymbol{\Lambda
}\mathbf{F}_0(0)-\mathbf{Q}\right) \left(\mathbf{p}^*_+(s)
\right)^{-1}\times\\ \times \left(s\mathbf{I}+\boldsymbol{\Lambda
}-\mathbf{Q} \right)^{-1}\int_{-\infty}^{0}\boldsymbol{\Pi}(dz)
\biggl(\left(\mathbf{I}-\mathbf{p}^*_+(s) \right)
    e^{\mathbf{R}^*_+(s)z}s
    \left(s\mathbf{I}+ \boldsymbol{\Lambda
    }\mathbf{F}_0(0)-\mathbf{Q}\right)^{-1}+\\
+\left(e^{i\alpha z}\left(\mathbf{C}-i\alpha
\mathbf{I}\right)-\left(\mathbf{I}-\mathbf{p}^*_+(s) \right)
e^{\mathbf{R}^*_+(s) z} \mathbf{C}\right)
\mathbf{p}^*_+(s)\left(\mathbf{R}^*_+(s)-i\alpha \mathbf{I}
   \right)^{-1}\boldsymbol{\Phi }^-(s,\alpha )
   \biggr).
\end{multline}
If $ m_1^0=\sum_{k=1}^m \pi_k (\lambda^1_k/c_k -\lambda ^2_k
m_k)>0$, then the following relations hold for $\eta
_{B,u}=\lim_{s\rightarrow 0}\eta _{B,u}(\theta_s)$
\begin{multline}\label{eq4.34}
\boldsymbol{\Phi }_{B,u}(\alpha )=\mathbf{E}e^{i\alpha \eta
  _{B,u}}=e^{i\alpha B}\left(\mathbf{I}-\mathbf{p}^*_+(0)
  \right)e^{-\mathbf{R}^*_+(0)v}\mathbf{p}^*_+
  \left(\boldsymbol{\Lambda }-\mathbf{Q} \right)^{-1}\int_{-\infty}^{0}\boldsymbol{\Pi }(dz)\times\\
\times\!\!\left(-e^{i\alpha
  z}\boldsymbol{\Psi  }^{-1}(\alpha )+\left(\mathbf{I}-\mathbf{p}^*_+(0)
  \right)e^{\mathbf{R}^*_+(0)z}\!\!\left(
  \left(\boldsymbol{\Lambda }\mathbf{F}_0(0)-\mathbf{Q} \right)^{-1}\!\!
  +\mathbf{C}\left(\mathbf{C}-i\alpha \mathbf{I} \right)^{-1}\boldsymbol{\Psi }^{-1}(\alpha )
  \right)\right),
\end{multline}
\begin{multline}\label{eq4.35}
\mathbf{P}\left\{\eta
_{B,u}=B\right\}=\left(\mathbf{I}-\mathbf{p}^*_+(0)
  \right)e^{-\mathbf{R}^*_+(0)v}\mathbf{p}^*_+
  \left(\boldsymbol{\Lambda }-\mathbf{Q} \right)^{-1}\times\\
\times\int_{-\infty}^{0}\boldsymbol{\Pi
}(dz)\left(\mathbf{I}-\mathbf{p}^*_+(0)
  \right)e^{\mathbf{R}^*_+(0)z}\left(
  \boldsymbol{\Lambda }\mathbf{F}_0(0)-\mathbf{Q} \right)^{-1}.
\end{multline}
\par For the dividend process $Y_{B,u}(\theta_s)$, the following
relations are true:
\begin{gather}
  \mathbf{E}e^{-\mu Y_{B,u}(\theta_s)}=\left(\mathbf{I}-\left(\mathbf{I}
  -\mathbf{p}^*_+(s) \right)e^{-\mathbf{R}^*_+(s)v}\mu
   \left(\mu \mathbf{I}+\mathbf{R}^*_+(s) \right)^{-1}
   \right)\mathbf{P}_s,\label{eq4.36}\\
   \mathbf{P}\left\{Y_{B,u}(\theta_s)=0\right\}=\mathbf{P}\left\{\xi
   ^+(\theta_s)<v\right\}=
   \mathbf{P}_s-\left(\mathbf{I}-\mathbf{p}^*_+(s) \right)
   e^{-\mathbf{R}^*_+(s)v}\mathbf{P}_s.\label{eq4.37}
\end{gather}
\end{theorem}
\noindent \textbf{Proof.} It follows from~\eqref{eq4.31} that the
characteristic function of $\eta _{B,u}(t)$ satisfies an integral
relation, which, after the Laplace--Karson transform, takes form
\begin{multline}\label{eq4.38}
    \boldsymbol{\Phi }_{B,u}(s,\alpha )=e^{i\alpha u}
    \mathbf{E}\left [e^{i\alpha \xi (\theta_s)},\xi ^+(\theta_s)<v
    \right]+\\
    +e^{i\alpha B}\mathbf{E}\left [  e^{-s T_1},T_1<\infty\right ]
     \left(\mathbf{I}-\mathbf{E}e^{-s\zeta _*}
     \right)\mathbf{P}_s
    +\mathbf{E}\left [  e^{-sT_1},T_1<\infty\right ]\mathbf{E}e^{-s\zeta _*}
    \widetilde{\boldsymbol{\Phi }}_B(s,\alpha ).
\end{multline}
\par According to~\cite[p.50]{Gusak2} and~\cite[p.43]{Gusak5}, we
get
\begin{gather}
  \mathbf{E}\left [  e^{-sT_1},T_1<\infty\right ]=\left(\mathbf{I}-\mathbf{p}^*_+(s) \right)
   e^{-\mathbf{R}^*_+(s)v},\notag\\
  \mathbf{E}e^{-s\zeta_* }=\left(s\mathbf{I}+\boldsymbol{\Lambda }\mathbf{F}_0(0)-\mathbf{Q} \right)^{-1}
    \boldsymbol{\Lambda }\mathbf{F}_0(0),\label{eq4.39}\\
  \mathbf{E}\left [  e^{i \alpha \xi (\theta_s)},\xi ^+(\theta_s)<v\right]
   =\mathbf{E}\left [  e^{i \alpha \xi^+ (\theta_s)},\xi
    ^+(\theta_s)<v\right]\mathbf{P}_s^{-1}\boldsymbol{\Phi }^-(s,\alpha)=\notag\\
   =\left(\!\left(\mathbf{C}-i\alpha \mathbf{I}\right)\!\mathbf{p}^*_+(s)-
   \left(\mathbf{I}-\mathbf{p}^*_+(s)\! \right)
    e^{\left(i\alpha \mathbf{I}-\mathbf{R}^*_+(s)
    \right)v}\mathbf{R}^*_+(s)\!\right)\!
   \left(\mathbf{R}^*_+(s)-i\alpha \mathbf{I}
   \right)^{-1}\!\boldsymbol{\Phi }^-(s,\alpha ).\notag
\end{gather}
Substituting these formulas into~\eqref{eq4.38}, we
obtain~\eqref{eq4.32} and the relation
\begin{multline}\label{eq4.40}
\widetilde{\boldsymbol{\Phi }}_B(s,\alpha )=e^{i\alpha
B}\left(\mathbf{I}-\int_{-\infty}^{0}d\mathbf{F}^1_0(z)\left(\mathbf{I}-\mathbf{p}^*_+(s)
\right)
e^{\mathbf{R}^*_+(s)z}\right.\times\\
\left.\times \left(s\mathbf{I}+\boldsymbol{\Lambda
}\mathbf{F}_0(0)-\mathbf{Q} \right)^{-1}\boldsymbol{\Lambda
}\mathbf{F}_0(0)\right)^{-1}
\left(\left(\int_{-\infty}^{0}e^{i\alpha
z}d\mathbf{F}^1_0(z)\left(\mathbf{C}-i\alpha
\mathbf{I}\right)\right.\right.-\\
\left.\left.-\int_{-\infty}^{0}d\mathbf{F}^1_0(z)\left(\mathbf{I}-\mathbf{p}^*_+(s)
\right) e^{\mathbf{R}^*_+(s) z} \mathbf{C}\right)\mathbf{p}^*_+(s)
\left(\mathbf{R}^*_+(s)-i\alpha \mathbf{I}
   \right)^{-1}\boldsymbol{\Phi }^-(s,\alpha )+\right.\\
\left.+\int_{-\infty}^{0}d\mathbf{F}^1_0(z)
   \left(\mathbf{I}-\mathbf{p}^*_+(s) \right)
    e^{\mathbf{R}^*_+(s)z}s
    \left(s\mathbf{I}+ \boldsymbol{\Lambda
    }\mathbf{F}_0(0)-\mathbf{Q}\right)^{-1}\right).
\end{multline}
The stochastic relations for $\tau ^+(x)$ (see~\cite[p.62]{Gusak2})
yield the following equation for $\mathbf{p}^*_+(s)$:
\begin{equation}\label{eq4.41}
  \left(s\mathbf{I}+\boldsymbol{\Lambda }-
  \mathbf{Q}\right)\left(\mathbf{I}-\mathbf{p}^*_+(s) \right)=
  \boldsymbol{\Lambda }\overline{\mathbf{F}}_0(0)+\boldsymbol{\Lambda
  }\mathbf{F}_0(0)\int_{-\infty}^{0}d\mathbf{F}^1_0(z)
  \left(\mathbf{I}-\mathbf{p}^*_+(s)\right)e^{\mathbf{C}\mathbf{p}^*_+(s)z}.
\end{equation}
Substituting~\eqref{eq4.41} into~\eqref{eq4.40}, we
receive~\eqref{eq4.33}.
Note that, for $m_1^0>0$: $\xi ^+ $ has a degenerate distribution.
Therefore, the first term in relation~\eqref{eq4.38} tends to $0$ as
$s\rightarrow 0$. Since $\left(\mathbf{I}-\mathbf{E}e^{-s\zeta _*}
\right)=s\left(s\mathbf{I}+\boldsymbol{\Lambda
    }\mathbf{F}_0(0)-\mathbf{Q} \right)^{-1}\mathbf{P}^{-1}_s$, the
    second term in~\eqref{eq4.38} tends to $0$ too. For the third term,
    by virtue of Theorem 3 in~\cite{Gusak5} and the first relation
    in~\eqref{eq4.2} we get
\begin{gather*}
  \lim_{s\rightarrow 0}s\left(\mathbf{p}^*_+(s) \right)^{-1}=\mathbf{p}^*_+,\\
\begin{split}
  \lim_{s\rightarrow 0}s^{-1}\mathbf{p}^*_+(s)\left(\mathbf{R}^*_+(s)
   -i\alpha \mathbf{I} \right)^{-1}\boldsymbol{\Phi }^-(s,\alpha
   )&=\lim_{s\rightarrow 0}(\mathbf{C}-i\alpha \mathbf{I})^{-1}
   \left(s\mathbf{I} -\boldsymbol{\Psi }(\alpha )\right)^{-1}=\\
   &=-(\mathbf{C}-i\alpha \mathbf{I})^{-1}\boldsymbol{\Psi}^{-1}(\alpha ).
\end{split}
\end{gather*}
Thus, passing to the limit as $s\rightarrow 0 $ in~\eqref{eq4.38} we
obtain~\eqref{eq4.34}. Relation~\eqref{eq4.36} and~\eqref{eq4.37}
follow from the representation: $
  Y^{kr}_{B,u}(t)\doteq\max\left(0, \xi ^+_{kr}(t)-v \right).
$

\end{document}